\documentclass[11pt,a4paper]{amsart}
\usepackage{amsmath,latexsym,amsfonts,amssymb,amsthm}
\usepackage{geometry}
\usepackage[
    colorlinks=true,
    linkcolor=blue,
    citecolor=red,
    urlcolor=red
]{hyperref}
\usepackage{mathrsfs}
\usepackage{graphicx,color}
\usepackage{tikz}
\usepackage{tikz-cd}
\usepackage{bbm}
\usepackage{comment}
\usepackage[new]{old-arrows}

\newtheorem{prop}{Proposition}[section]
\newtheorem{thm}[prop]{Theorem}

\newtheorem{lem}[prop]{Lemma}

\theoremstyle{definition}

\newtheorem{rem}[prop]{\it Remark}

\newtheorem*{claim*}{Claim}

\newcommand{\bP}{\mathbb{P}}
\newcommand{\bC}{\mathbb{C}}

\newcommand{\Hom}{\mathrm{Hom}}

\newcommand{\an}{\mathrm{an}}
\newcommand{\Bl}{\mathrm{Bl}}

\newcommand{\Z}{\mathbb Z}

\newcommand{\C}{\mathbb C}

\begin{document}

\title[Isolated symplectic singularities admitting projective resolutions]{Isolated symplectic singularities admitting symplectic resolutions}
\author{Duo Li}
\address{Duo Li, School of Mathematics, Sun-Yat Sen University,
Zhuhai, Guangdong 519082, PR China}
\email{liduo5@mail.sysu.edu.cn}

\author{Yanjie Li}
\address{Yanjie Li, Academy of Mathematics and Systems Science, Chinese Academy of Sciences,
55 Zhong Guan Cun East Road, Beijing 100190, China,
University of Chinese Academy of Sciences, Beijing, China}
\email{liyanjie241@mails.ucas.ac.cn}

\begin{abstract}
We give a proof that an isolated algebraic symplectic singularity of dimension $2n\ge 4$
admitting a projective symplectic resolution is analytically isomorphic to the minimal
nilpotent orbit closure in $\mathfrak{sl}_{n+1}$ at the origin. 
\end{abstract}

\maketitle

\section{Introduction}

Symplectic singularities, introduced by Beauville \cite{Beauville2000},
are normal singularities whose smooth locus carries a holomorphic
symplectic form extending to a regular two-form on a resolution. They arise naturally not only in algebraic geometry, but also in geometric representation theory and mathematical physics.
Symplectic singularities may admit a particularly interesting
type of resolution, called a \emph{symplectic resolution}: a resolution $\pi:Y \rightarrow X$ of a symplectic singularity $X$ is called a symplectic resolution if the pull back of the symplectic form on $X_{\mathrm{reg}}$ by $\pi$ can be extended to a symplectic form on $Y$. 

The classification of symplectic singularities is a natural and
challenging problem. In dimension two, symplectic singularities are precisely the rational double points, and their minimal
resolutions are symplectic \cite[\S2]{Beauville2000}.
In higher dimensions, the picture is considerably more complicated. Among symplectic singularities, those with an isolated singular
point form a natural first class to study. Even this class,
however, contains a rich variety of examples.
Besides the classical examples of isolated symplectic quotient
singularities and minimal nilpotent orbit closures
\cite[\S2]{Beauville2000}, several new constructions have emerged
in recent years.
Bellamy, Bonnaf\'e, Fu, Juteau, Levy and Sommers constructed an
infinite family of four-dimensional isolated symplectic
singularities with trivial local fundamental group
\cite{BellamyEtAl2023}.
Namikawa subsequently obtained such examples in every even
dimension at least four by toric hyperK\"ahler reduction
\cite{Namikawa2025Isolated}.
Further examples continue to appear; see, for instance,
\cite{Berry2025}.

By contrast, up to local analytic isomorphism, the only known
isolated symplectic singularities of dimension at least four
admitting projective symplectic resolutions belong to a single
series: the closures of the minimal nonzero nilpotent orbits
of type $A$,
\[
    \overline{\mathcal O}_{\min}(\mathfrak{sl}_{n+1})
    =\{A\in\mathfrak{sl}_{n+1}\mid
        \operatorname{rk}A\leq1,\ A^2=0\},
    \qquad n\geq2.
\]
The variety $\overline{\mathcal O}_{\min}(\mathfrak{sl}_{n+1})$ has dimension $2n$, is singular only at the
origin, and admits the Springer resolution
$T^*\mathbb P^n\longrightarrow \overline{\mathcal O}_{\min}(\mathfrak{sl}_{n+1})$,
which contracts the zero section to the origin; see
\cite[\S2.3]{Fu2006}.
It is expected that these examples exhaust all isolated
complex algebraic symplectic singularities of dimension at least
four admitting projective symplectic resolutions.

In dimension four, this classification was established by Wierzba
and Wi\'sniewski \cite[Theorem~1.1]{WierzbaWisniewski2003}.
In higher dimensions, Cho, Miyaoka and Shepherd-Barron also studied this problem and stated
the corresponding classification for contractions of projective
symplectic manifolds
\cite[Theorem~8.3]{ChoMiyaokaShepherdBarron2002}.
Brion and Fu subsequently proved the classification under the
assumption that the singularity admits a contracting
$\mathbb C^*$-action of symplectic weight one
\cite[Proposition~6]{BrionFu2015}. Combined with the recent work \cite{XuZhuang2026, Odaka2026}, we extend the method of \cite{BrionFu2015} to obtain a proof the following theorem.

\begin{thm}\label{thm:main}
Let $(X,o)$ be an isolated algebraic symplectic singularity of dimension $2n$, and
suppose that it admits a projective symplectic resolution. If $n\ge2$, then we have an isomorphism of analytic germs
\[
(X,o)_{\an}\simeq
\bigl(\overline{\mathcal O}_{\min}(\mathfrak{sl}_{n+1}),0\bigr)_{\an}.
\]
The analytic germ of the resolution along this fibre is isomorphic to the germ of
\[
T^*\bP^n\longrightarrow \overline{\mathcal O}_{\min}(\mathfrak{sl}_{n+1})
\]
along its zero section.
\end{thm}

Our proof combines recent results on the conical structure of
symplectic singularities with the approach of Brion
and Fu in \cite{BrionFu2015}. Using Odaka's hyperk\"ahler cone theorem
\cite[Corollary~1.3]{Odaka2026}, we first show that every isolated
symplectic singularity admits an analytic model with an effective
good $\mathbb C^*$-action of symplectic weight one or two. The weight-one case follows from
\cite[Proposition~6]{BrionFu2015}.
When the weight is $2$, following the strategy in \cite{WierzbaWisniewski2003} and \cite{ChoMiyaokaShepherdBarron2002}, we first identify the normalizations of the
exceptional components with $\mathbb P^n$.
The main step is to show that the equivariant normalization
$\nu:\mathbb P^n\to D$ of a component containing the source of the action is an isomorphism, which is Proposition \ref{prop:smooth-component}.

\subsection*{AI disclosure} 
We discussed with ChatGPT many times. The proofs obtained with the assistance of AI were reconstructed and revised by the authors. We take full responsibility for the correctness and content of the paper.

\subsection*{Acknowledgements} 
We are deeply grateful to Professor Baohua Fu for introducing us to this problem, for providing the key ideas underlying all the proofs in this paper, and for his continued support and numerous fruitful discussions throughout this project.

\section{$\mathbb{C}^*$-action of weight $1$ or $2$}
In this section, we apply the result in \cite{Odaka2026} to show that any isolated symplectic singularity is locally isomorphic to a symplectic singularity with an effective good $\mathbb{C}^*$-action of weight $1$ or $2$. Recall that a $\bC^*$-action on $X=\operatorname{Spec}R$ is called
\emph{good} if $R=\bigoplus_{j\ge0}R_j$ with $R_0=\bC$,
and \emph{effective} if its kernel is trivial, see \cite[Introduction]{Namikawa2013}.

\begin{prop}\label{thm:weight 1 or 2}
  Let $(X,x)$ be an isolated symplectic singularity. Then the analytic germ of $(X,x)$ is analytically isomorphic to that of $(X_0,x_0)$, where $X_0$
is a symplectic variety endowed with an effective good
algebraic $\C^*$-action and an algebraic symplectic form $\omega_0$
satisfying
\(
    \lambda^*\omega_0=\lambda^\ell\omega_0
\)
with $\ell\in\{1,2\}$.
\end{prop}
\begin{proof}
By \cite[Corollary~1.3]{Odaka2026}, the germ $(X,x)$
is analytically isomorphic to the germ of the vertex of an affine hyperk\"ahler
cone $X_0$. Let $r$ be the radius function on $X_0$.  Since the singularity is isolated, the link $L=\{r=1\}$ is smooth
and compact.

Let $I,J,K$ denote the hyperk\"ahler complex structures on $X_0$.
The vector fields $I(r\partial_r)$, $J(r\partial_r)$ and
$K(r\partial_r)$ are tangent to $L$ and generate a Lie algebra isomorphic to $\mathrm{sp}(1)$.
Since $L$ is compact, every real linear combination of these
vector fields is complete. By the Lie--Palais theorem,
the infinitesimal action integrates to an isometric
$\mathrm{Sp}(1)$-action on $L$.
In particular, the Reeb field $\xi=I(r\partial_r)$ satisfies
$\exp(2\pi\xi)=\mathrm{id}$.

For a torus eigenfunction $f$ with $r\partial_r f=af$, we have
$\xi f=\sqrt{-1} af$, so $\exp(2\pi \xi) = \mathrm{id}$ implies $a\in\Z$.
Then we have an integral grading
\( \C[X_0]=\bigoplus_{j\in\mathbb Z}R_j.\)
The grading associated with the cone structure is positive, we have
$R_j=0$ for $j<0$ and $R_0=\C$.
Thus this grading defines a good algebraic
$\C^*$-action on $X_0$.
By \cite[Propositions~4.1 and~4.2]{Odaka2026},
$X_0$ admits an algebraic symplectic form $\omega_0$
parallel with respect to the hyperk\"ahler cone metric.
Let $\nabla$ be the the Levi-Civita connection of the cone metric. Since $\nabla R=\mathrm{Id}$, we have
$\mathcal L_R\omega_0=2\omega_0$.
Thus the induced algebraic $\C^*$-action satisfies
$\lambda^*\omega_0=\lambda^2\omega_0$.
Consequently, the kernel of this action is contained in
$\{\pm1\}$. Passing to the effective quotient gives weight
$\ell=2$ or $\ell=1$, respectively.
\end{proof}

\begin{rem}\label{rem:resolution}
    Suppose that $(X,o)_{\an}\simeq(X_0,x_0)_{\an}$
and that $X$ admits a projective symplectic resolution.
Transporting the resolution germ to $(X_0,x_0)$ and
gluing it with $X_0^{\an}\setminus\{x_0\}$ gives a
proper analytic crepant resolution of $X_0^{\an}$.

This resolution on $X_0^{\an}$ is projective. Indeed, since $X_0^{\an}$
is Stein and contractible, the exponential sequence gives
$H^2(X_0^{\an},\mathcal O^*)=0$.
By the Mayer--Vietoris sequence for the cover
$X_0^{\an}=U\cup(X_0^{\an}\setminus\{x_0\})$,
where $U$ is a sufficiently small neighbourhood of $x_0$,
a relatively ample line bundle on $f^{-1}(U)$ extends,
after twisting by the pullback of a line bundle on $U$,
to a globally defined $f$-ample line bundle.
Thus a relatively ample line bundle near $x_0$
extends after twisting by a line bundle from the base. By \cite[Proposition~A.7]{Namikawa2008},
with the correction in
\cite[Lemma~3.1 and the following Remark]{Namikawa2026},
it algebraizes to a $\mathbb C^*$-equivariant
projective crepant resolution of $X_0$.
Since $X_0$ is symplectic, this resolution is
symplectic by \cite[Proposition~1.6]{Fu2006}.
\end{rem}

\section{Normalizations of the exceptional components}

Let $X$ be an isolated affine symplectic variety of dimension $2n$ with singular locus $\{o\}$, and let $\omega \in H^0(X,\Omega_X^{[2]})$ be the reflexive symplectic $2$-form. Assume that $X$ carries a good $\bC^*$-action of symplectic weight $\ell>0$.
Let $\pi:Y\to X$ be a projective symplectic resolution, put
$\Omega=\pi^*\omega$. Let $E=\pi^{-1}(o)_{\mathrm{red}}$ be the exceptional fiber.
The action lifts to $Y$ by \cite[Proposition~A.7]{Namikawa2008}. We keep this notation through Section~4.

We first recall the description of the normalization of the components of $E$ proved in \cite{ChoMiyaokaShepherdBarron2002}. The following lemma is \cite[Corollary 8.7]{ChoMiyaokaShepherdBarron2002} without the assumption that $Y$ is projective.

\begin{lem}\label{lem:quotient}
Let $q:V\to B$ be a proper bimeromorphic map from a
holomorphic symplectic manifold of dimension $2n$ to a
normal Stein space, and suppose that $q$ is an isomorphism
over $B\setminus\{o\}$.
Let $D$ be a positive-dimensional irreducible component of
$q^{-1}(o)_{\mathrm{red}}$, with normalization $\nu:P\to D$.
Assume that $P$ is projective.
Then $P\simeq\bP^n/\Gamma$ for a finite subgroup
$\Gamma\subset\mathrm{PGL}_{n+1}(\bC)$ whose nonidentity
elements have no divisorial fixed locus.
\end{lem}

\begin{proof}
We follow the proof of
\cite[Corollary~8.7]{ChoMiyaokaShepherdBarron2002},
using Ran's deformation estimate in the noncompact
analytic setting.
Put $F=q^{-1}(o)_{\mathrm{red}}$.
As $V$ is symplectic, the morphism $q$ is crepant and $B$ has canonical singularities.
Thus $E$ is rationally chain connected by
\cite[Corollary~1.2(i) and Remark~1.5]{Fujino2026RCC},
and $P$ is uniruled.
Moreover, $F$ is isotropic by
\cite[Proposition~2.17]{AmerikVerbitsky2024}, so $\dim P =  n$.

Fix an ample line bundle on $P$ and choose a dominating irreducible component $H$ of rational curves of minimal degree.
Let $\overline H\subset\operatorname{Chow}(P)$ be the closure
of its image, and let $\overline H_x$ denote the closed
subfamily of cycles passing through $x\in P$.
By minimality, $\overline H_x$ consists of irreducible
reduced rational curves for general $x$.

Choose a general parametrization $f:\bP^1\to P$ of a member
of $H$, lying on no other component of the corresponding
fixed-degree Hom space, such that $\nu f(\bP^1)$ meets
$D_{\mathrm{reg}}$ away from the other components of $F$.
By \cite[Corollary~5.1]{Ran1995},
\[
\dim_{[\nu\circ f]}\operatorname{Hom}(\bP^1,V)\geq2n+1.
\]
Since $B$ is Stein and $q$ is an isomorphism off $F$, nearby nonconstant deformations of $\nu\circ f$ remain in $F$, and hence in $D$ by our choice of $f$.
After normalizing their parameter space, they lift
uniquely to $P$.
Their degree is unchanged, so the generality of $f$
ensures that the lifted curves remain in $H$.
By considering the $\mathrm{PGL_2}(\mathbb{C})$-action, we obtain $\dim H\geq2n-2$.

For general $x\in P$,
\cite[Proposition~2.4(2)]{ChoMiyaokaShepherdBarron2002}
therefore gives
\[
2n-2\leq\dim H
\leq\dim P+\dim\operatorname{Locus}(\overline H_x)-2
\leq2\dim P-2\leq2n-2.
\]
Hence $\operatorname{Locus}(\overline H_x)=P$.
By properness of the two-point evaluation map,
$\overline H$ is doubly dominant.
It is maximal by the choice of $H$ as an irreducible
component of $\operatorname{RatCurves}^n(P)$.
The conclusion now follows from
\cite[Theorem~0.2]{ChoMiyaokaShepherdBarron2002}.
\end{proof}

\section{Classification of the exceptional fiber for $l=2$}
In this section, we prove that the exceptional fiber $E = \pi^{-1}(o)_{\mathrm{red}}$ is isomorphic to $\mathbb{P}^n$ when the $\mathbb{C}^*$-action on $X$ is of weight $l=2$. We first analyze the $\mathbb{C}^*$-equivariant normalization map $\nu: \mathbb{P}^n \rightarrow D$ of some component $D$ containing the component of $Y^{\mathbb{C}^*}$ which is the \emph{source} of $Y$ with respect to the $\mathbb{C}^*$-action. Then following the stragey used in \cite{ChoMiyaokaShepherdBarron2002} and \cite{WierzbaWisniewski2003}, we prove that if $E$ has a smooth component, then it is isomorphic to $\mathbb{P}^n$.

\subsection{Bia\l{}ynicki-Birula decomposition}
Since the action on $X$ is contracting and $\pi$ is proper, the limit $\lim_{t\to0}t\cdot y$ exists for every $y\in Y$, and
$Y^{\bC^*}\subset E$ is projective. Thus $Y$ is semiprojective, and we use its
Bia\l ynicki--Birula decomposition as in \cite[Section~1]{HauselRodriguezVillegas2015}.
For a connected component $F\subset Y^{\bC^*}$, we denote the \emph{repelling cell} by
\[
Y_F^-=\bigl\{y\in Y\mid\lim_{t\to\infty}t\cdot y\in F\bigr\}.
\]
This is a smooth locally closed subvariety, an affine-space bundle over $F$.
Since the action on $X$ is good and $E$ is projective, we have $E=\bigcup_F Y_F^-$.
For $y\in F$, let
\[
T_yY=\bigoplus_{a\in\mathbb Z}(T_yY)_a
\]
be the weight-space decomposition of the $\C^*$-representation on $T_yY$.
By Homogeneity, we have $$\Omega((T_yY)_a,(T_yY)_b)=0$$ unless $a+b=\ell$.
Since $\ell>0$, the tangent space $T_yY_F^-=\bigoplus_{a\le0}(T_yY)_a$ is isotropic.
Hence
\begin{equation}\label{eq:dimension-bound}
\dim Y_F^-=\sum_{a \leq 0}\dim (T_yY)_a\le n.
\end{equation}

\begin{prop}\label{prop:normalization}
Every irreducible component of $E$ has normalization $\bP^n$.
\end{prop}

\begin{proof}
Let $D\subset E$ be an irreducible component. By Lemma~\ref{lem:quotient}, its normalization is $P=\bP^n/\Gamma$; write $q:\bP^n\to P$ for the quotient map.
Since $E$ is a finite union of the cells $Y_F^-$, some cell $U=Y_F^-$ is dense
in $D$. Now we fix the notation $F$ to denote the corresponding fixed component of this dense cell. Then $U = Y_F^{-}$ is smooth and locally closed, hence open in $D$, and we identify it with its inverse image in $P$.

The composite rational map $\bP^n\dashrightarrow U\longrightarrow F$
is dominant. Thus the smooth projective variety $F$ is unirational, so
$\pi_1(F)=1$ by \cite{Serre1959}.
As $U\to F$ is a bundle of affine space, we have $\pi_1(U)=\pi_1(F)=1$.
The finite map $q^{-1}(U)\to U$ is \'etale in codimension one, and hence
everywhere by purity of the branch locus. Since $q^{-1}(U)$ is a nonempty open
subset of $\bP^n$, it is connected. This cover therefore has degree one,
so $\Gamma=1$.
\end{proof}

Let $M\subset Y^{\bC^*}$ be the \emph{source}, namely
the fixed component whose attracting cell for $t\to0$ is open in $Y$.
Then $M$ is a smooth projective variety and we have $(T_pY)_a = 0$ for $a<0$ and $p \in M$. Set $m=\dim M$.

Assume that $\ell=2$, we consider the following weight decomposition
\begin{equation}\label{eq:source-weights}
T_Y|_M=E_0\oplus E_1\oplus E_2,\qquad
E_0=T_M,\quad E_2\simeq\Omega_M^1,
\end{equation}
where the last isomorphism is induced by the symplectic paring. The bundle $E_1$ is symplectic. 

\begin{prop}\label{prop:smooth-component}
If the weight $l=2$, then some irreducible component of $E$ containing $M$ is smooth and isomorphic to
$\bP^n$.
\end{prop}

\begin{proof}
Choose a component $D$ of $E$ containing $M$, and let $\nu:\bP^n\to D$
be its normalization (see Proposition \ref{prop:normalization}). 
The $\mathbb{C}^*$-action lifts uniquely to the normalization
$\nu:P=\bP^n\to D$.
We identify $P=\bP(V)$ with
\[
V=\bigoplus_{a\ge0}V_a,\qquad
t\cdot v=t^av\quad(v\in V_a),
\qquad V_0\ne0,
\]
where the integral weights are taken with respect to the
original $\bC^*$-action and normalized so that the smallest
weight is zero.

Since $M$ is the source, its repelling cell is $Y_M^-=M$.
As $\nu$ is finite and equivariant, $\nu^{-1}(M)$ consists
of fixed points. If a fixed component $\bP(V_a)$ with $a>0$
mapped into $M$, the repelling cell of $\bP(V_a)$ would
contain nonfixed points mapped into $Y_M^-=M$, a contradiction.
Thus $\nu^{-1}(M)=\bP(V_0)$, which is the source of $P$. We have a finite surjective morphism
\( h:=\nu|_{\mathbb{P}(V_0)}:\mathbb{P}(V_0)\simeq\bP^m\longrightarrow M\).

Now we prove that if $V_a \ne 0$, then $a=0,1$. If there is a weight $a>2$ such that $V_a \ne 0$, every point of a line joining $V_a$ to
$V_0$ would be fixed by $\mu_a \subseteq \mathbb{C}^*$. Its image would lie in the connected component $S_a$ of $Y^{\mu_a}$ containing $M$. As $a>2$ and the weights of $T_Y|_M$ contains only $0,1,2$ (see (\ref{eq:source-weights})), we have $T_yS_a = (T_yY)^{\mu_a} = (T_y Y)_0 = T_yM$ for $y \in M$. Hence $S_a=M$, contrary to $\nu^{-1}(M)=\mathbb{P}(V_0)$. We next exclude the case $V_2 \ne 0$. Let $s = \dim V_2$ and let $S$ be the connected component of $Y^{\mu_2}$ containing $M$. The subspace $\mathbb{P}(V_0 \oplus V_2) \simeq \bP^{m+s}$ maps into the connected component $S$ of $Y^{\mu_2}$ containing
$M$. For $y \in M$, we have $(T_yY)_0 = T_yM$ and $T_yS = (T_y Y)^{\mu_2} = (T_yY)_0 \oplus (T_yY)_2 $. Since $\Omega$ has weight $2$, the variety $S$ is a symplectic manifold and therefore $\dim S =2\dim (T_yY)_0 = 2m$. On the other hand, we have $\nu^*\Omega=0$ as $H^0(\mathbb{P}^n,\Omega_{\mathbb{P}^n}^2)=0$. The image $\nu(\bP^{m+s})$ is therefore isotropic in $S$, so $m+s\le m$. Thus $s=0$.

If $\dim M = n$, then $D=M$ is smooth and isomorphic to $\bP^n$. Suppose that $n>m$ and we set $r:=n-m$. We denote two fixed subspaces in $P = \mathbb{P}(V)$ by
\[
F_{-}:=\mathbb{P}(V_0)=\bP^m,\qquad F_{+}:=\mathbb{P}(V_1)=\bP^{r-1}.
\]
Let $F\subset Y^{\bC^*}$ be the connected component containing $\nu(F_+)$. Since $\nu$ is $\mathbb{C}^*$-equivariant, we have
\(D\setminus M=\nu(\bP^n\setminus F_-)\subset Y_F^-.\)
The cell $Y_F^-$ is irreducible. Since $\dim(D\setminus M)=n$,
\eqref{eq:dimension-bound} shows that $D\setminus M$ is dense in it.
As $D$ is closed, $Y_F^-\subset D$. Since $F\ne M$ (as $\nu^{-1}(M)=\mathbb{P}(V_0)$), we have
$Y_F^-\cap M=F \cap M =\varnothing$, so $Y_F^-=D\setminus M$.
This is a smooth variety, and the normalization map $\nu$ gives
\begin{equation}\label{eq:normalization-isom-off-source}
\nu:\bP^n\setminus F_-\xrightarrow{\ \sim\ }D\setminus M.
\end{equation}

To prove that $\nu$ is an isomorphism, we analyze its behavior in the normal directions to $F_{-}$. The other open subset $\bP^n\setminus F_+$ is the total space of $N=N_{F_-/\bP^n}\simeq\mathcal O_{\bP^m}(1)^{\oplus r}$, with fiberwise $\mathbb{C}^*$-action of weight one. Recall the decomposition (\ref{eq:source-weights}), then we denote the weight-one part of $d\nu$ along the zero section $F_{-}$ of $ \mathrm{Tot}(N) \simeq \mathbb{P}^n \backslash F_+$ by
\begin{equation}\label{eq:normal-differential}
a:N\longrightarrow h^*E_1.
\end{equation}
We claim that $a$ is fiberwise injective.
Suppose that $a_x(v)=0$ for some $x\in F_-$ and
$0\ne v\in N_x$. Consider the equivariant morphisms
\[
    \gamma_\pm:\mathbb A^1\longrightarrow Y,\qquad
    \gamma_\pm(t)=\nu(x,\pm tv).
\]
They both map $0$ to $h(x)$, and
$d\gamma_{\pm,0}(\partial_t)=\pm a_x(v)=0$. Suppose that their formal germs at $0$ differ, and let $k$
be the first order at which they differ.
Since $\gamma_-(t)=\gamma_+(-t)$, the difference
$\gamma_+^*u-\gamma_-^*u$ contains only odd powers of $t$
for every $u\in\mathcal O_{Y,h(x)}$.
Its linear term vanishes, so $k$ is odd and $k\geq3$.
Writing
$\gamma_+^*u-\gamma_-^*u
=t^k\delta_k(u)\pmod{t^{k+1}}$,
the leading difference $\delta_k$ is a derivation at $h(x)$,
and hence defines a nonzero vector in $T_{h(x)}Y$. As $\nu$ is $\mathbb{C}^*$-equivariant, we have $\delta_k\in(T_{h(x)}Y)_k$, contradicting \eqref{eq:source-weights} as $k\geq3$.

Therefore $\gamma_+$ and $\gamma_-$ agree formally at $0$ and hence agree on a Zariski neighbourhood of $0$ by \cite[Tag~00IQ]{stacks-project}.
For any nonzero $t$ in this neighbourhood, we obtain
$\nu(x,tv)=\nu(x,-tv)$.
However, $(x,tv)$ and $(x,-tv)$ are distinct points of
$\operatorname{Tot}(N)\setminus F_-\subset P\setminus F_-$,
contradicting \eqref{eq:normalization-isom-off-source}.
Hence $a$ is a subbundle inclusion.

We now aim to show that if $h: F_{-} = \mathbb{P}(V_0) \rightarrow M$ is of degree larger than $1$, then there is two  distinct points $x_1,x_2 \in F_{-}$ such that there exists non zero elements in the fibers $v_1 \in N_{x_1}, v_2 \in N_{x_2}$ such that $(x_1,v_1)$ and $(x_2,v_2)$ are mapped to the same point by $\nu$ (we view $(x_i,v_i) (i=1,2)$ as points in $\mathrm{Tot}(N) \cong \mathbb{P}^n \backslash F_{+}$), contradicting the isomorphism (\ref{eq:normalization-isom-off-source}). We need the following property of equivariant morphisms that to compare the branches over $M$.

\begin{lem}\label{lem:weight-two-rigidity}
Let $q:B\to M$ be a dominant morphism from an integral variety
$B$, and let $N_B$ be a vector bundle on $B$, with the
fiberwise $\bC^*$-action of weight one and the trivial action
on $B$. Suppose that two equivariant morphisms
$f_1,f_2:\operatorname{Tot}(N_B)\longrightarrow Y$
satisfy $f_1|_B=f_2|_B=q$, where $B$ is the zero section
and $q$ is viewed as a morphism to $Y$.
If their normal differentials agree, i.e.
\[
    (df_1)|_{N_B}=(df_2)|_{N_B}:N_B\longrightarrow q^*E_1,
\]
and assume that $\operatorname{Hom}_{\mathcal O_B}\bigl(\operatorname{Sym}^2N_B,q^*E_2\bigr)=0$, then $f_1=f_2$.
\end{lem}

\begin{proof}
Let $\mathcal I$ be the ideal of the zero section in
$\operatorname{Tot}(N_B)$.
The morphisms $f_1,f_2$ agree modulo $\mathcal I^2$.
Suppose that they differ on the formal completion along $B$,
and let $k\geq2$ be the first order at which they differ.
Thus they agree modulo $\mathcal I^k$ but not modulo
$\mathcal I^{k+1}$.

The difference of the pullbacks of local regular
functions on $Y$ by $f_1$ and $f_2$ modulo $\mathcal I^{k+1}$ defines a derivation relative to $q$ with values in $\mathcal I^k/\mathcal I^{k+1}$. Then we have a nonzero equivariant homomorphism
\[
    q^*\Omega_Y^1
    \longrightarrow
    \mathcal I^k/\mathcal I^{k+1}
    \simeq\operatorname{Sym}^kN_B^\vee.
\]
Dually, this defines a nonzero equivariant homomorphism
\[
    \delta_k:\operatorname{Sym}^kN_B
    \longrightarrow q^*(T_Y|_M).
\]
The source has weight $k$, while the target has only weights
$0,1,2$ by \eqref{eq:source-weights}.
Hence $k=2$ and the image of $\delta_2$ lies in $q^*E_2$,
contradicting the assumption $\operatorname{Hom}_{\mathcal O_B}\bigl(\operatorname{Sym}^2N_B,q^*E_2\bigr)=0$.

Thus $f_1$ and $f_2$ agree formally along $B$.
By \cite[Tag~00IQ]{stacks-project}, we have $f_1=f_2$ on a neighbourhood of $B$.
Using the fiberwise $\bC^*$-action, we obtain
$f_1=f_2$ on $\operatorname{Tot}(N_B)$.
\end{proof}

We return to the proof of Proposition \ref{prop:smooth-component}. Suppose $\deg h >1$ and $m>0$ (if $m=0$, then $\mathbb{P}^0 \rightarrow M$ is an isomorphism). Let $T \subset (F_{-} \times_M F_{-})_{\mathrm{red}}$ be an irreducible component dominating $M$ that is not contained in the diagonal. We denote by $p_i: T \rightarrow F_{-}(i=1,2)$ the corresponding projection. Let $(x_1,x_2) \in T$ be a point in $T$ and we set $z=h(x_1)=h(x_2)$. We construct an isomorphism $\varphi: p_1^*N \xrightarrow{\sim} p_2^*N$ as follows. Set $q:=h p_1 = hp_2$. Since $a: N \rightarrow h^*E_1$ is an injection of vector bundles, the pull-back $a_i:=p_i^*a: p_i^*N \rightarrow q^*E_1 (i=1,2)$ is also an injection of vector bundles. As the image of $a$ is a Lagrangian subbunlde of $h^*E_1$, we have the following exact sequence
$$0 \rightarrow p_2^*N \xrightarrow{a_2} q^*E_1 \xrightarrow{} (p_2^*N)^{\vee} \rightarrow 0.$$
From $N = \mathcal{O}_{\mathbb{P}^m}(1)^{\oplus r}$ and that the morphisms $p_i(i=1,2)$ are finite, we have $\Hom(p_1^*N,(p_2^*N)^{\vee})$ = 0. Then the injection $p_1^*N \rightarrow q^*E_1$ of vector bundles induces an isomorphism $\varphi: p_1^*N \xrightarrow{\sim} p_2^*N$ such that $p_2^*a \circ \varphi = p_1^*a$.

To apply Lemma \ref{lem:weight-two-rigidity} for $B:=T$, $q: B \rightarrow M$ and $N_B = p_1^*N$, we verify the vanishing $\Hom_{\mathcal{O}_B}(\mathrm{Sym}^2(p_1^*N),q^*E_2) = 0$. As $h: \mathbb{P}^m \rightarrow M$ is finite, the differential map $h^*\Omega_M^1 \rightarrow \Omega_{\mathbb{P}^m}^1$ is an isomorphism at generic point of $\mathbb{P}^m$. Pull back by the dominant morphism $p_1$, we also get a morphism $q^*\Omega_M^1 \rightarrow p_1^*\Omega_{\mathbb{P}^m}^1$ which is an isomorphism over the generic point. Since $E_2 \simeq \Omega^1_M$ by the duality induced by the symplectic form (\ref{eq:source-weights}), we therefore have an injection of sheaves $q^*E_2 \hookrightarrow p_1^*\Omega_{\mathbb{P}^m}^1 \hookrightarrow p_1^*\mathcal{O}_{\mathbb{P}^m}(-1)^{\oplus m+1}$, where the second injection follows from the Euler sequence $0 \rightarrow \Omega^1_{\mathbb{P}^m} \rightarrow \mathcal{O}_{\mathbb{P}^m}(-1)^{\oplus m+1} \rightarrow \mathcal{O}_{\mathbb{P}^m} \rightarrow 0$. If $\Hom_{\mathcal{O}_B}(\mathrm{Sym}^2(p_1^*N),q^*E_2) \ne 0$, we will get a non-zero element in $\Hom((p_1^*\mathcal{O}_{\mathbb{P}^m}(1))^{\oplus 2},p_1^*\mathcal{O}_{\mathbb{P}^m}(-1))$, which is impossible as $p_1$ is finite.

Now we apply the lemma \ref{lem:weight-two-rigidity}. A point in the total space $\mathrm{Tot}(N_B)$ is a pair $((x_1,x_2),v)$, where $(x_1,x_2)$ is a point in $B$ and $v$ is a vector in $N_{x_1}$. We now define the morphisms $f_1,f_2: \mathrm{Tot}(N_B) \rightarrow Y$ by 
\[ f_1((x_1,x_2),v)= \nu(x_1,v),~ f_2((x_1,x_2),v)=\nu(x_2,\varphi(v))\]
where we view $(x_1,v)$ and $(x_2,\varphi(v))$ as points in $\mathrm{Tot}(N)\simeq \mathbb{P}^n\backslash F_{+}$. Then 
\[df_1|_{p_1^*N} = p_1^*a = p_2^*a \circ \varphi = df_2|_{p_1^*N}\]
and we apply the lemma \ref{lem:weight-two-rigidity} to get $f_1 = f_2$. Choose $(x_1,x_2)\in T$ off the diagonal and
$0\ne v\in N_{x_1}$.
Then $(x_1,v)$ and $(x_2,\varphi(v))$ are distinct points
of $P\setminus F_-$ with the same image under $\nu$,
contradicting \eqref{eq:normalization-isom-off-source}.
Thus $\deg h=1$.
Since $h$ is finite and $M$ is smooth, $h$ is an isomorphism.

Now $\nu: P \rightarrow D$ is a bijection, we check that the differential map $d\nu_x$ along $F_{-}$ is injective. At $x \in F_{-}$, we consider the weight decomposition of the tangent space $T_xP = T_x F_{-} \oplus N_x$. Let $(u,v) \in T_xP$ be a tangent vector at $x$, where $u \in T_xF_{-}$ is of weight zero and $v \in  N_x$ is of weight $1$. Then $d\nu_x(u,v) =( dh_x(u),a_x(v)) \in T_{h(x)}M \oplus (E_1)_{h(x)}$ by the $\mathbb{C}^*$-equivariance of $\nu$. As $h$ is an isomorphism and $a$ is an injection, the map $d\nu_x$ is injective. Then $d\nu$ is everywhere injective and $\nu$ is therefore an isomorphism.
\end{proof}

\subsection{Irreducibility of the exceptional fibre}

In this subsection $\ell>0$ is arbitrary and $n\ge2$. We prove the following proposition by adapting the Mukai-flop
argument in \cite[proof of Lemma~8.11]{ChoMiyaokaShepherdBarron2002},
also used in dimension four in \cite[proof of Proposition~3.2]{WierzbaWisniewski2003}.
\begin{prop}\label{prop:irreducibility}
If a component $P$ of $E$ is smooth and isomorphic to $\bP^n$,
then $E=P$.
\end{prop}

We first show that $P$ meets every other component in finitely many points.

\begin{lem}\label{lem:finite-intersection}
Let $P\simeq\bP^n$ be a smooth subvariety of a smooth complex variety $Z$,
with $N_{P/Z}\simeq\Omega_P^1$. If $\nu:F=\bP^n\to Z$ is finite onto its
image and $\nu(F)\not\subset P$, then $\nu^{-1}(P)$ is finite.
\end{lem}

\begin{proof}
Since $N_{P/Z}^{\vee}\simeq T_P$ is ample, Grauert's
contraction theorem \cite[Section~3, Satz~8]{Grauert1962}
gives an analytic neighbourhood $U\subset Z^{\an}$ of $P$
and a proper holomorphic map $\rho:U\to V$ to a Stein space,
such that $\rho(P)=\{p\}$ and
$\rho:U\setminus P\to V\setminus\{p\}$ is an isomorphism.
Every compact irreducible curve in $U$ is contained in $P$.

Put $S=\nu^{-1}(P)$ and suppose that it contains an
irreducible curve $C$.
By compactness of $C$, we have $\nu(g(C))\subset U$
for every $g\in\mathrm{PGL}_{n+1}$ sufficiently close
to the identity in the analytic topology.
Since $\nu$ is finite, $\nu(g(C))$ is a compact
irreducible curve, and hence lies in $P$.
Thus $g(C)\subset S$ for all such $g$.

Fix $x\in C$.
As $g$ varies near the identity, the points $g(x)$
fill an analytic neighbourhood of $x$ in $F$.
Hence $S$ contains a nonempty analytic open subset of $F$,
so $S=F$, contrary to $\nu(F)\not\subset P$.
Thus $S$ contains no curve, and hence is finite.
\end{proof}

We recall the following Lemma proved in \cite[Lemma 8.2]{ChoMiyaokaShepherdBarron2002} showing that the Mukai flop of $P$ exists in the analytic category.
\begin{lem}[{\cite[Lemma 8.2]{ChoMiyaokaShepherdBarron2002}}]\label{lem:mukai-flop}
Let $P\simeq\bP^n$ be a smooth component of $E$, with $n\geq2$,
and let $P^*$ denote the dual projective space parametrizing
hyperplanes in $P$.
Let $b:(\Bl_PY)^{\an}\to Y^{\an}$ be the blow-up along $P$.
Its exceptional divisor $Q=b^{-1}(P)$ is naturally identified
with the incidence variety
\[
Q\simeq\{(x,H)\in P\times P^*:x\in H\},
\]
under which $b|_Q$ is the first projection. There exist a smooth holomorphic symplectic manifold $Y^+$
containing $P^*$ and a proper bimeromorphic map $b^+$
fitting into
\[
Y^{\an}\xleftarrow{\ b\ }(\Bl_PY)^{\an}
\xrightarrow{\ b^+\ }Y^+,
\]
such that $b^+|_Q$ is the second projection
$(x,H)\mapsto H$ and $b^+$ induces an isomorphism
between $(\Bl_PY)^{\an}\setminus Q$ and $Y^+\setminus P^*$. Moreover, $\pi^{\an}\circ b$ descends to a proper
bimeromorphic map $\pi^+:Y^+\to X^{\an}$, which iss an isomorphism over $X^{\an}\setminus\{o\}$.
\end{lem}

\begin{proof}[Proof of Proposition~\ref{prop:irreducibility}]
Suppose that $E\ne P$. As $E$ is connected, we can choose another component $D(\ne P)$ of $E$ such that $D \cap P$ is non-empty. Let $\nu: F \simeq \mathbb{P}^n \rightarrow D$ be the normalization map of $D$. As $P$ is a Lagrangian subvariety in $Y$, we have $N_{P/Y} \simeq \Omega_P^1$. By Lemma \ref{lem:finite-intersection}, $S = \nu^{-1}(P)$ is an non-empty finite subset.

We consider the Mukai flop of $P$ as in Lemma \ref{lem:mukai-flop} and use the notations there. Let $\widetilde{D} \subset \Bl_P Y$ be the strict transform of $D$ in $(\Bl_P Y)^{\an}$ and set $D^+ = b^{+}(\widetilde{D})$. Set $\mathcal{I} = \nu^{-1}\mathcal{I}_P \cdot \mathcal{O}_F$ and consider the normalized blow up $g: \widehat{F}=(\Bl_{\mathcal{I}}F)^{\nu} \rightarrow F$. We prove that $\widehat{F} \rightarrow D^+$ is the noramlization map. As $F \rightarrow Y$ is finite, the morphism $\Bl_{\mathcal{I}}F \rightarrow \Bl_PY$ is also finite. We note the image of the latter map is $\widetilde{D}$. Then the morphism $\widehat{F} \rightarrow \widetilde{D}$ is finite and birational (it is an isomorphism over $b^{-1}(D_{\mathrm{reg}} \backslash P)$). The map $\widetilde{D} \rightarrow D^+$ is also finite and birational: for $ H\in P^* \subset Y^{+}$, the fiber of $\widetilde{D} \rightarrow D^+$ over $H$
$$\widetilde{D} \cap (b^+)^{-1}(H) \subseteq \{(x,H) \in Q~|~x\in D \cap P\},$$
which is finite as $D \cap P $ is finite. Then we have $\widehat{F} \simeq (D^+)^{\nu}$. As $D^{+}$ is an irreducible component of $\pi^{-1}(o)_{\mathrm{red}}$. Applying Lemma \ref{lem:quotient} to $\pi^+$ and $D^+$, we have $\widehat{F} \simeq \mathbb{P}^n/\Gamma$ for some finite subgroup $\Gamma \subset \mathrm{PGL}_{n+1}(\mathbb{C})$ and hence $\rho(\widehat{F})=1$. Then $g: \widehat{F}=(\Bl_{\mathcal{I}}F)^{\nu} \rightarrow F$ is a finite map and hence an isomorphism. It follows that $\mathcal{I} = \nu^{-1}\mathcal{I}_P \cdot \mathcal{O}_F$ is an invertible sheaf with $\mathrm{supp}(\mathcal{O}_F / \mathcal{I})= S=\nu^{-1}(P)$, which is impossible since $S$ is finite and $\dim F = n \geq 2$.
\end{proof}

We now prove our main theorem.

\begin{proof}[Proof of Theorem~\ref{thm:main}]
By Proposition \ref{thm:weight 1 or 2}, we may assume that the isolated symplectic singularity $X$ admits an effective good $\mathbb{C}^*$-action of weight $1$ or $2$, and that the symplectic resolution $\pi: Y \rightarrow X$ is equivariant with respect to this $\mathbb{C}^*$-action (by \cite[Proposition~A.7]{Namikawa2008}, see Remark \ref{rem:resolution}). If the weight of the $\mathbb{C}^*$-action is $1$, we can apply \cite[Proposition 6]{BrionFu2015}. If the weight of the $\mathbb{C}^*$-action is $2$, the classification follows from Proposition \ref{prop:smooth-component}, Proposition \ref{prop:irreducibility} and \cite[Lemma 8.2]{ChoMiyaokaShepherdBarron2002} (see also \cite[Theorem 1.2]{AmerikVerbitsky2024}).
\end{proof}

\bibliographystyle{alpha}
\bibliography{isolated_symplectic_classification}

\end{document}